\documentclass[twocolumn]{IEEEtran}

\def\BibTeX{{\rm B\kern-.05em{\sc i\kern-.025em b}\kern-.08em
    T\kern-.1667em\lower.7ex\hbox{E}\kern-.125emX}}

\usepackage{latexsym}
\usepackage{amssymb}

\def\qed{\hfill$\Box$}

\newcommand{\bydef}{\mbox{$ \;\stackrel{\triangle}{=}\; $}}

\newcommand{\geqd}{\mbox{$ \;\stackrel{(d)}{\geq}\; $}}
\newcommand{\geqe}{\mbox{$ \;\stackrel{(e)}{\geq}\; $}}
\newcommand{\eqa}{\mbox{$ \;\stackrel{(a)}{=}\; $}}
\newcommand{\eqb}{\mbox{$ \;\stackrel{(b)}{=}\; $}}
\newcommand{\eqc}{\mbox{$ \;\stackrel{(c)}{=}\; $}}

\newcommand{\RL}{{\mathbb R}}

\newcommand{\VAR}{\mbox{\rm Var}} 
 
\newcommand{\signs}{\mbox{\scriptsize sign}}
\newcommand{\Ahat}{\mbox{$\hat{A}$}}
\newcommand{\Ahatn}{\mbox{$\hat{A}^n$}}

\newcommand{\Dmax}{\mbox{$D_{\rm max}$}}

\newcommand{\DminP}{\mbox{$D_{\rm min}^{P,Q}$}}
\newcommand{\DmaxP}{\mbox{$D_{\rm max}^{P,Q}$}}

\newcommand{\LA}{\mbox{$\Lambda$}}

\newcommand{\la}{\mbox{$\lambda$}}

\def\be{\begin{eqnarray}}
\def\ee{\end{eqnarray}}
\def\ben{\begin{eqnarray*}}
\def\een{\end{eqnarray*}}

\input{epsf}

\newtheorem{theorem}{Theorem}
\newtheorem{proposition}{Proposition}
\newtheorem{lemma}{Lemma}
\newtheorem{corollary}{Corollary}
\newtheorem{example}{Example}
\begin{document}

\title{Critical Behavior in Lossy Source Coding}

\author{Amir Dembo, Ioannis Kontoyiannis\thanks{A. Dembo
is with the Department of
        Mathematics and the Department of
        Statistics,
        Stanford University,
        Stanford, CA 94305.
        Email: amir@stat.stanford.edu.
I. Kontoyiannis is with the
Department
        of Statistics, Purdue University,
        1399 Mathematical Sciences Building,
        W.~Lafayette, IN 47907-1399.
        Email:  yiannis@stat.purdue.edu.}
\thanks{A.D.'s research was
supported in part by
NSF grant \#DMS-0072331.
I.K.'s research was
supported in part by
NSF grant \#0073378-CCR
and by a grant from the Purdue
Research Foundation.}}

\markboth{IEEE Transactions On Information Theory, Vol. XX, No. Y, Month
2000}
{Dembo and Kontoyiannis: Critical Redundancy in Lossy Source Coding} 

\maketitle

\begin{abstract}
The following critical phenomenon
was recently discovered.
When a memoryless source is
compressed using a variable-length
fixed-distortion code, the fastest
convergence rate of the (pointwise)
compression ratio to $R(D)$ is
either $O(\sqrt{n})$ or $O(\log n)$.
We show it is always $O(\sqrt{n})$,
except for discrete, uniformly
distributed sources.
\end{abstract}


\section{Introduction}
\PARstart{S}{uppose} 
that data is produced by a stationary memoryless
source $\{X_n\;;\;n\geq 1\}$, so that the $X_i$
are independent and identically distributed
(IID) random variables with common distribution
$P$. We will assume throughout that the $X_i$
take values in the {\em source alphabet} $A$,
where $A$ is a subset of $\RL$,
and that the {\em reproduction alphabet}
$\Ahat$ is a finite subset of $\RL$,
say $\Ahat=\{a_1,a_2,\ldots,a_k\}$.

The main objective of data compression 
is to find efficient approximate representations 
for data $x_1^n=(x_1,x_2,\ldots,x_n)$ 
generated from 
the source $X_1^n=(X_1,X_2,\ldots,X_n)$. 
Specifically, we wish to represent 
each source string $x_1^n$ by a corresponding 
string $y_1^n=(y_1,y_2,\ldots,y_n)$ 
taking values in the reproduction alphabet
$\Ahat$, so that the distortion between each
$x_1^n$ and its representation lies within 
some fixed allowable range. For our purposes, 
distortion is measured by a family of 
single-letter distortion measures,
\be
\rho_n(x_1^n,y_1^n)=\frac{1}{n}
\sum_{i=1}^n\rho(x_i,y_i)
\;\;\;\;\;\;\;\;\;x_1^n\in A^n,\;\;y_1^n\in\Ahatn,
\label{eq:singlelet}
\ee
where $\rho:A\!\times\!\Ahat\to[0,\infty)$ 
is a fixed nonnegative function.

We consider 
{\em variable-length block codes operating 
at a fixed distortion level}, that is, codes 
$C_n$ defined by triplets 
$(B_n,\phi_n,\psi_n)$ where:

\newcounter{abc}
\begin{list}{($\alph{abc})$}
  {\usecounter{abc}
	\setlength{\itemindent}{0cm}
	\setlength{\topsep}{0.2cm}
	\setlength{\labelwidth}{0.3cm}
	\setlength{\leftmargin}{0.4in}
	\setlength{\labelsep}{0.2cm}
	\setlength{\rightmargin}{0.1in}
	\setlength{\parsep}{0.5ex plus0.2ex minus0.1ex}
	\setlength{\itemsep}{0ex plus0.2ex}
  }
	\item $B_n$ is a subset of $\Ahatn$ called the {\em codebook};
	\item $\phi_n:A^n\to B_n$ is a map called the {\em encoder};
	\item $\,\psi_n:B_n\to \{0,1\}^*$ is a prefix-free
		representation of the elements of $B_n$ by
		finite-length binary strings.
\end{list}
For a fixed distortion level $D\geq 0$, 
the code $C_n=(B_n,\phi_n,\psi_n)$ is said to 
{\em operate at distortion level $D$} 
\cite{kieffer:91} 
if it encodes 
each source string with distortion $D$ 
or less:
$$\rho_n(x_1^n,\phi_n(x_1^n))\leq D
\;\;\;\;\mbox{for all}\;\;x_1^n\in A^n.$$
Our main quantity of interest here is the 
description length of a block code $C_n$,
expressed by its length function 
\mbox{$\ell_n:A^n\to{\mathbb N}$:}
$$\ell_n(x_1^n)=\;
\mbox{length of}\;[\psi_n(\phi_n(x_1^n))].$$
Broadly speaking, the smaller the 
description length, the better the code.

Shannon's celebrated source coding 
theorem 
states that, for an 
arbitrary sequence of block codes
$\{C_n=(B_n,\phi_n,\psi_n)\;;\;
n\geq 1\}$
operating at distortion level $D$,
the expected compression 
ratio $E[\ell_n(X_1^n)]/n$ is
asymptotically bounded below by 
the rate-distortion function
$R(D)$:
\ben
\liminf_{n\to\infty}\,
\frac{E[\ell_n(X_1^n)]}{n}\geq R(D)\;\;\;\;
\mbox{bits per symbol.}
\een
Moreover, Shannon showed that there exist
codes achieving the above lower bound with 
equality; see Shannon's 1959 paper \cite{shannon:59} 
or Berger's classic text 
\cite{berger:book}.
A stronger version of Shannon's theorem
was proved by Kieffer in 1991 
\cite{kieffer:91}, where it is
shown that 
the rate-distortion function is a
{\em pointwise}
asymptotic lower bound for $\ell_n(X_1^n)$:
\be
\liminf_{n\to\infty}\,
\frac{\ell_n(X_1^n)}{n}\geq R(D)\;\;\;\;
\mbox{with prob. 1.}
\label{eq:kieffer}
\ee
In \cite{kieffer:91} it is also
demonstrated that the bound in 
(\ref{eq:kieffer}) can be achieved 
with equality.

The following
refinement to Kieffer's result was
recently given in \cite{kontoyiannis-red:00}:

\vspace{0.1in}

\begin{quote}
(POINTWISE REDUNDANCY): {\em For any sequence
of block codes $\{C_n\}$ with associated
length functions $\{\ell_n\}$,
operating at distortion level $D$, 
\be
\ell_n(X_1^n) 
	&\geq& nR(D)+\sum_{i=1}^nf(X_i)-2\log n
	\nonumber\\
	& & \hspace{0.4in}
	\mbox{eventually, with prob. 1,}
\label{eq:LBD}
\ee
where $f:A\to\RL$ is a bounded function 
depending on $P$ and $D$ but {\em not} 
on the codes $\{C_n\}$, such that $E_P[f(X_1)]=0$. 
Moreover, there exist codes $\{C_n,\ell_n\}$ 
that achieve
\be
\ell_n(X_1^n)
	&\leq& nR(D)+\sum_{i=1}^nf(X_i)+5\log n
	\nonumber\\
	& & \hspace{0.4in}
	\mbox{eventually, with prob. 1.}
\label{eq:UBD}
\ee
}
\end{quote}

\noindent
[{\sl cf.} Theorems~4 and~5 and eq.~(18)
in \cite{kontoyiannis-red:00}; 
above and throughout the paper,
`$\log$' denotes the logarithm taken to base 2
and `$\log_e$' denotes the natural logarithm.]
The function $f$ is defined precisely in Section~III;
here we just mention the following interpretation: 
If we write $\tilde{f}(x)=f(x)+R(D)$, 
then $\tilde{f}$ can be expressed in a natural 
way in terms of familiar information theoretic 
quantities. In particular, 
$E(\tilde{f}(X_1))=R(D)$,
its variance $\sigma^2=\VAR(\tilde{f}(X_1))$
is the ``minimal coding variance'' of the 
source with distribution $P$
\cite{kontoyiannis-red:00}, and in 
the case of lossless compression 
(as $\,D\downarrow 0\,$), $\tilde{f}(x)$
reduces to $-\log P(x)$.

The above result says that, for any source distribution
$P$ and any sequence of codes $\{C_n\}$ operating 
at distortion level $D$, the ``pointwise 
redundancy'' in the description lengths 
of the codes $C_n$, namely, the difference 
between $\ell_n(X_1^n)$ and 
the optimum $nR(D)$ bits,
is essentially bounded below by the sum 
of the IID, bounded, zero-mean random 
variables $f(X_i)$. So there
are two possibilities:
\begin{list}{$\bullet$}
{
\setlength{\itemindent}{0cm}
        \setlength{\topsep}{0.2cm}
        \setlength{\labelwidth}{0.3cm}
        \setlength{\leftmargin}{0.4in}
        \setlength{\labelsep}{0.2cm}
        \setlength{\rightmargin}{0.1in}
        \setlength{\parsep}{0.5ex plus0.2ex minus0.1ex}
        \setlength{\itemsep}{0ex plus0.2ex}
}
\item
Either the random variables $f(X_i)$ are
non-constant, in which case the best
achievable pointwise redundancy rate will
be of order $O(\sqrt{n})\,$
(by the central limit theorem and the
upper and lower bounds in
(\ref{eq:LBD}) and (\ref{eq:UBD}));
\item
or the random variables $f(X_i)$ are
equal to zero with probability one,
in which case the best achievable
pointwise redundancy is no more than
$(5\log n)$ bits, eventually
(by (\ref{eq:UBD})).
\end{list}
To be more precise, in the first case when 
the random variables $f(X_i)$ are not constant, 
the central limit theorem implies that the term
$\sum_{i=1}^nf(X_i)$ is of order $O(\sqrt{n})$ 
in probability, and therefore, by (\ref{eq:LBD}) 
and (\ref{eq:UBD}), the best achievable pointwise
redundancy will also be of order $O(\sqrt{n})$ 
in probability. [In a similar fashion, the law of 
the iterated logarithm implies that the 
{\em pointwise} fluctuations of the best achievable 
pointwise redundancy will be of order 
$O(\sqrt{n\log_e \log_e n})$; see 
\cite[Section~I]{kontoyiannis-red:00} 
for a more detailed discussion. Also
the contrast between the pointwise and
the expected redundancy rate is interpreted 
and commented on in 
\cite[Remark~3, p.139]{kontoyiannis-red:00}.]

Our purpose in this paper is to characterize
exactly when each one of the above two cases
occurs, namely, when the minimal pointwise 
redundancy is $O(\sqrt{n})$ and when it is 
$O(\log n)$. In the next section we show that 
it is almost never the case that $f(X_1)=0$ 
with probability one, so the minimal pointwise 
redundancy is typically of order $\sqrt{n}$.
In particular, in the common case 
when the $X_i$ take values in 
a finite alphabet $A=\Ahat$, then 
(under mild conditions) we show that 
$f(X_1)=0$ with probability one if 
and only if the $X_i$ are uniformly 
distributed.

Before stating our main results 
(Theorems~1, 2 and~3 in the next section)
in detail, we recall the following
representative examples from 
\cite{kontoyiannis-97} and
\cite{kontoyiannis-red:00}.

\vspace{0.1in}

\begin{example}[Lossless Compression]
For a source $\{X_n\}$ with
distribution $P$ on the finite 
alphabet $A$,
a lossless 
code $C_n$ is a prefix-free
map $\psi_n:A^n\to\{0,1\}^*$.
[Or, to be pedantic, in our setting
a lossless code is a code operating 
at distortion level $D=0$ with 
respect to Hamming distortion.]
In this case
the function $f$ 
has the simple form
\be
f(x)=-\log P(x)-H(P)
\label{eq:lossless}
\ee
where $H(P)=E_P[-\log P(X_1)]$ is
the entropy of $P$,
and 
the lower bound (\ref{eq:LBD}) is
simply
\be
    \ell_n(X_1^n)
	&\geq& 
	nH(P)+\sum_{i=1}^nf(X_i)-2\log n
	\nonumber\\
	&=& 
	-\log P(X_1^n) -2\log n
	\label{eq:barron}\\
	& & 
	\hspace{0.6in}
	\mbox{eventually, with prob. 1.}
	\nonumber
\ee
The lower bound (\ref{eq:barron}) 
is a well-known
information-theoretic fact
called Barron's lemma (see
\cite{algoet:thesis}\cite{barron:thesis}
and the discussion in \cite{kontoyiannis-red:00}).
It says that the description lengths
$\ell_n(X_1^n)$ of an arbitrary sequence
of codes are (eventually
with probability 1) bounded
below by the idealized Shannon code lengths
$-\log P(X_1^n)$, up to terms of order $\log n$.
From (\ref{eq:lossless})
it is obvious that $f(X_1)=0$ with probability 
one if and only if $P$ is the uniform
distribution on $A$. 
\end{example}

\vspace{0.1in}

\begin{example}[Binary Source, 
Hamming Distortion]$\;$\\
This is
the simplest non-trivial 
lossy example.
Suppose $\{X_n\}$ is a binary source
with Bernoulli($p$) 
distribution for some $p\in(0,1/2]$. 
Let $A=\Ahat=\{0,1\}$
and take $\rho$ to be the Hamming
distortion measure, $\rho(x,y)=0$ 
when $x=y$, and equal to 1 otherwise.
For each fixed $D\in(0,p)$ it is shown 
in \cite{kontoyiannis-red:00} that
$$f(x)=-\log\left(\frac{P(x)}{1-D}\right)
	-E_P\left[
		-\log\left(\frac{P(X_1)}{1-D}\right)
		\right],$$
from which it is again obvious that
$f(X_1)=0$ with probability one
if and only if $p=1/2$, i.e., if and
only if $P$ is the uniform distribution 
on $A=\{0,1\}$.
\end{example}

\vspace{0.1in}

In a third example presented in 
\cite{kontoyiannis-red:00} 
it is also found that $f(X_1)=0$ 
with probability one if and only 
if $P$ is the uniform distribution,
and the natural question is raised as to
whether this pattern persists in general.
In the next section we answer this 
question by showing (in Theorem~1 and Corollary~1)
that for a source distribution $P$
on a finite alphabet, $f(X_1)$ can be 
equal to zero with probability one for 
at most finitely many distortion levels 
$D$, {\em unless} $P$ is the uniform 
distribution and $\rho$ is a ``permutation'' 
distortion measure.
In Theorems~2 and~3 and in Corollary~2 
the continuous case is considered, 
and it is shown that when $P$ is a 
continuous distribution it essentially
never happens that $f(X_1)=0$ with 
probability one.
Section~III contains the proofs of 
Theorems~1, 2 and~3 and Corollaries~1 and~2.

\section{Results}

Suppose that the source alphabet
$A$ is an arbitrary (Borel) subset of $\RL$,
and let $P$ be a (Borel) probability measure 
on $\RL$, supported on $A$ (the special cases 
when $P$ is purely discrete or purely continuous
are considered separately below). Let 
$\Ahat=\{a_1,a_2,\ldots,a_k\}$ be the
finite reproduction alphabet of size $k$.
Given an arbitrary, bounded,
nonnegative function
$\rho:A\!\times\!\Ahat\to[0,M]$ (for some
finite constant $M$),
define a sequence of single-letter
distortion measures 
$\rho_n:A^n\!\times\!\Ahatn\to[0,M]$
as in (\ref{eq:singlelet}). 
Throughout the paper,
we make the usual assumption:
\be
\sup_{x\in A}\min_{y\in\Ahat}\rho(x,y)=0.
\label{eq:basic}
\ee
[See, e.g.,
\cite[p.26]{berger:book} or
\cite[Ch.13,~ex.4]{cover:book};
if (\ref{eq:basic}) is
not satisfied, for example when $A$ is
an interval of real numbers, $\hat{A}$
is a finite set, and $\rho(x,y)=(x-y)^2$,
we may consider the
distortion measure
$\rho'(x,y)=\rho(x,y)-\min_{z\in\hat{A}}\rho(x,z)$
instead.]
For $D\geq 0$,
the {\em rate-distortion function} 
of a memoryless source with
distribution $P$ is
\be
R(D) = \inf_{(X,Y)} I(X;Y)
\label{eq:rd-def}
\ee
where the infimum is over all
jointly distributed random
variables $(X,Y)$ with
values in $A\!\times\!\Ahat$
such that $X$ has distribution $P$
and $E[\rho(X,Y)]\leq D$;
$I(X;Y)$ denotes the mutual information
(in bits) between $X$ and $Y$ 
(see \cite{berger:book} for 
more details).
Under our assumptions,
the rate-distortion function $R(D)$ 
is a convex, nonincreasing function
of $D\geq 0$, and it is finite for all $D$.

For a fixed distribution $P$ on $A$,
let
\ben
\Dmax \;=\;\Dmax(P)\;=\; \min_{y\in\hat{A}}\; E_{P}[\rho(X,y)]
\een
and recall that
$R(D)=0$ for $D\geq\Dmax$
(see, e.g., Proposition~1 in
Section~III).
In order to avoid the trivial 
case when $R(D)$ is identically
zero we assume that $\Dmax>0$,
and from now on we restrict our
attention to the interesting range 
of distortion levels $D\in(0,\Dmax)$.

\subsection{The Discrete Case: $A=\Ahat$}

We first consider the most 
common case where the source
$\{X_n\}$ takes values in 
a finite alphabet
$A=\Ahat=\{a_1,a_2,\ldots,a_k\}$.
Suppose that $\{X_n\}$ are IID with 
common distribution $P$ on $A$,
and assume, without loss of generality,
that $P_i=P(a_i)>0$ for all $i=1,\ldots,k$.
Given a distortion measure $\rho$,
write $\rho_{ij}$ for $\rho(a_i,a_j)$.
We assume throughout this section
that $\rho$ is symmetric,
i.e., that $\rho_{ij}=\rho_{ji}$ 
for all $i,j$,
and also that
$\rho_{ij}=0$ if and only if
$i=j$. We call $\rho$ a {\em permutation
distortion measure}, if all
rows of the matrix $(\rho_{ij})_{i,j=1,\ldots,k}$
are permutations of one another
(which, by symmetry, is equivalent 
to saying that all columns
are permutations of one another).

Recall that the minimal pointwise redundancy
is of order $O(\log n)$ if and only if $f(X_1)=0$
with probability one; otherwise it is
$O(\sqrt{n})$. 
Our first result says that 
the rate cannot be $O(\log n)$ for 
many distortion levels $D$,
unless the distribution $P$ is uniform
in which case the rate is 
$O(\log n)$ for all
distortion levels $D$.

\vspace{0.1in}

\begin{theorem}$\;$\\
\indent
$(a)\;$ If $P$ is the uniform distribution 
on $A$ and $\rho$ is a 
permutation distortion measure,
then $f(X_1)=0$ with probability one
for all $D\in(0,\Dmax)$.\\
\indent
$(b)\;$
If $f(X_1)=0$ with probability one for
a sequence of distortion values 
$D_n\in(0,\Dmax)$ such that $D_n\downarrow 0$,
then $P$ is the uniform distribution and
$\rho$ is a permutation distortion
measure, and therefore $f(X_1)=0$ with probability one
for {\em all} $D\in(0,\Dmax)$.
\end{theorem}

\vspace{0.1in}

As we mentioned above, the rate-distortion
function $R(D)$ is convex for $D\in(0,\Dmax)$.
If it is {\em strictly} convex
(as it is usually the case -- see the
discussion in \cite[Chapter~2]{berger:book}),
then Theorem~1 can be strengthened to 
the following.

\vspace{0.1in}

\begin{corollary}
Suppose $R(D)$ is strictly convex 
over the range $D\in(0,\Dmax)$. If $f(X_1)=0$ 
with probability one for infinitely 
many $D\in(0,\Dmax)$ then $P$ is the 
uniform distribution and $\rho$ is 
a permutation distortion measure, 
and therefore $f(X_1)=0$ with 
probability one
for {\em all} $D\in(0,\Dmax)$.
\end{corollary}

\vspace{0.1in}

{\sl Remark.} In the examples presented
in the previous section it turned out
that either $f(X_1)=0$ with probability
one for {\em all} $D$, or it was never
the case. But it may happen that $f(X_1)=0$ 
with probability one only for a few isolated 
values of $D$, while $P$ is {\em not} the
uniform distribution. Such an example is
given after Lemma~3 in Section~III-B.

\subsection{The Continuous Case: $A=\RL$}

Here we take $A=\RL$ and we
assume that the distribution $P$ of
the source has a positive density $g$ 
(with respect to Lebesgue measure), 
or, more generally, that there exists
a (nonempty) open interval $I\subset\RL$ 
on which $P$ has an absolutely continuous
component with density $g$ such that
$g(x)>0$ for $x\in I$.
Since the
reproduction alphabet 
$\Ahat=\{a_1,a_2,\ldots,a_k\}$
is finite, 
given a distortion measure $\rho$
we can write
$r_j(x)=\rho(x,a_j)$
for 
all $1\leq j\leq k$ and all
$x\in A.$
We assume that for all $j$ the 
functions $r_j$ are continuous on $I$.
For convenience we also define,
for $j=0$, $r_j(x)\equiv 0$ on $I$.

Our next result gives a sufficient
condition on the distortion measures
$r_j$, under which the best 
redundancy rate in (\ref{eq:kieffer}) 
can {\em never} be $O(\log n)$.

\vspace{0.1in}

\begin{theorem}
If for every $\la<0$
the functions 
$$e^{\lambda r_j(\cdot)},\;\;\;\;\;j=0,1,\ldots,k$$
are linearly independent on $I$, then
$f(X_1)$ cannot be equal
to zero with probability one for any
distortion level $D\in(0,\Dmax).$
\end{theorem}

\vspace{0.1in}

Next we provide a somewhat
simpler set of conditions, 
under which we get a weaker 
conclusion. Theorem~3 says
that the best redundancy 
rate in (\ref{eq:kieffer}) cannot
be $O(\log n)$ for many distortion
levels $D$.

\vspace{0.1in}

\begin{theorem}
Under either one of the following two
conditions,
$f(X_1)$ cannot be equal
to zero with probability one for distortion 
levels $D>0$ arbitrarily close to zero.\\
\indent
($a$)
There exist (distinct) points $\{x_0,x_1,\ldots,x_k\}$
in $I$ such that, for all $0\leq i\neq j\leq k$, with
$j\neq 0$, we have $r_j(x_j) > r_j(x_i).$\\
\indent
($b$)
There exist (distinct) points $\{x_0,x_1,\ldots,x_k\}$
in $I$ such that, for every permutation $\pi$ of
the indices $\{0,1,\ldots,k\}$ with $\pi$
not equal to the identity, we have
$$\sum_{j=0}^k r_j(x_j)\neq
\sum_{j=0}^k r_j(x_{\pi(j)}).$$
\end{theorem}

\vspace{0.1in}

Although the conditions of Theorems~2 
and~3 may seem unusual, they are natural 
and generally easy to verify. To illustrate 
this, we present below two simple examples.

\vspace{0.1in}

\begin{example}[Mean-Squared Error]
Suppose $P$ has a positive density
on the interval $I=[-2,2]$, let
$\Ahat$ consist of the two reproduction
points $\pm 1$, and let $\rho$ be
the mean-squared error distortion
measure. Recall that, to satisfy 
(\ref{eq:basic}), $\rho(x,y)$ is 
actually defined by 
$$\rho(x,y)=(x-y)^2-\min\{(x-1)^2, (x+1)^2\}.$$
The corresponding distortion functions
$r_1(x)=\rho(x,-1)$ and
$r_2(x)=\rho(x,+1)$ are shown
in Figure~1. Here, condition~(a) 
of Theorem~3 is easily seen to hold with
$x_0=0,$ $x_1=2$ and $x_2=-2$.
\end{example}

\begin{figure}[ht]
\centerline{\epsfxsize 6cm \epsfbox{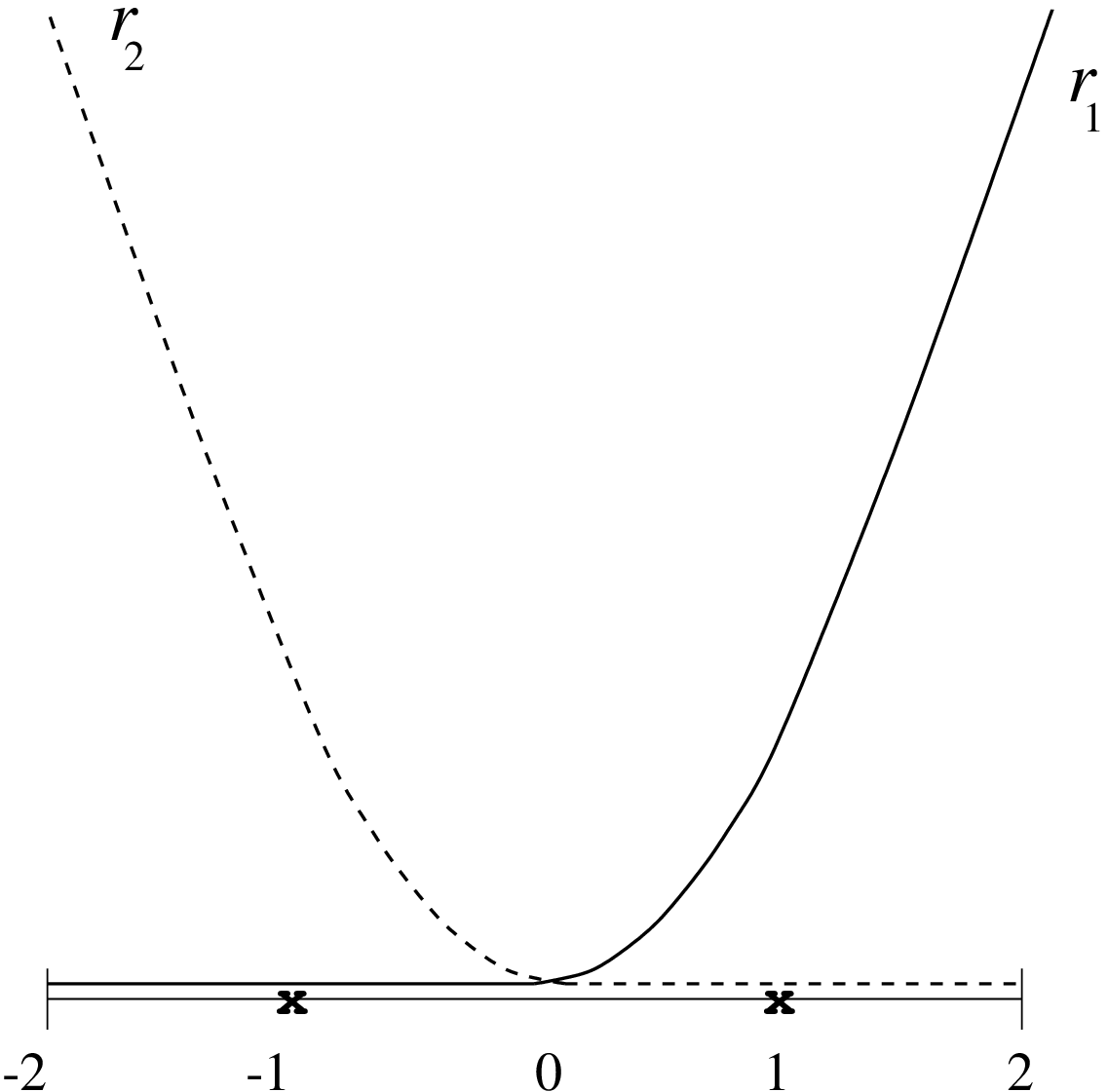}}
\caption{Distortion measure in Example~3. 
Reproduction points are shown as {\bf x}'s.}
\end{figure}

\begin{example}[$L^1$ Distance]
Suppose $P$ has a positive density
on the interval $I=[0,6]$, let
$\Ahat=\{1,3,5\},$ 
and take $\rho$ to be 
the normalized $L^1$ distance
$|x-y|$ adjusted so that 
(\ref{eq:basic}) is satisfied;
the resulting functions $r_j(\cdot)$
are shown in Figure~2. Here it is
easy to verify that the 
condition of Theorem~2 is satisfied, 
i.e., that the functions 
$\{e^{\lambda r_j(\cdot)}\;;\;0\leq j\leq 3\}$
are linearly independent on $I$. For this
it suffices to observe that $e^{\lambda r_1}$
and $e^{\lambda r_3}$ are linearly independent
on $[2,4]$ (essentially because the functions
$e^{\lambda x}$ and $e^{-\lambda x}$ are
linearly independent on $[0,2]$),
and that $e^{\lambda r_2}$ is not constant
outside $[2,4]$.
\end{example}

\begin{figure}[ht]
\centerline{\epsfxsize 7cm \epsfbox{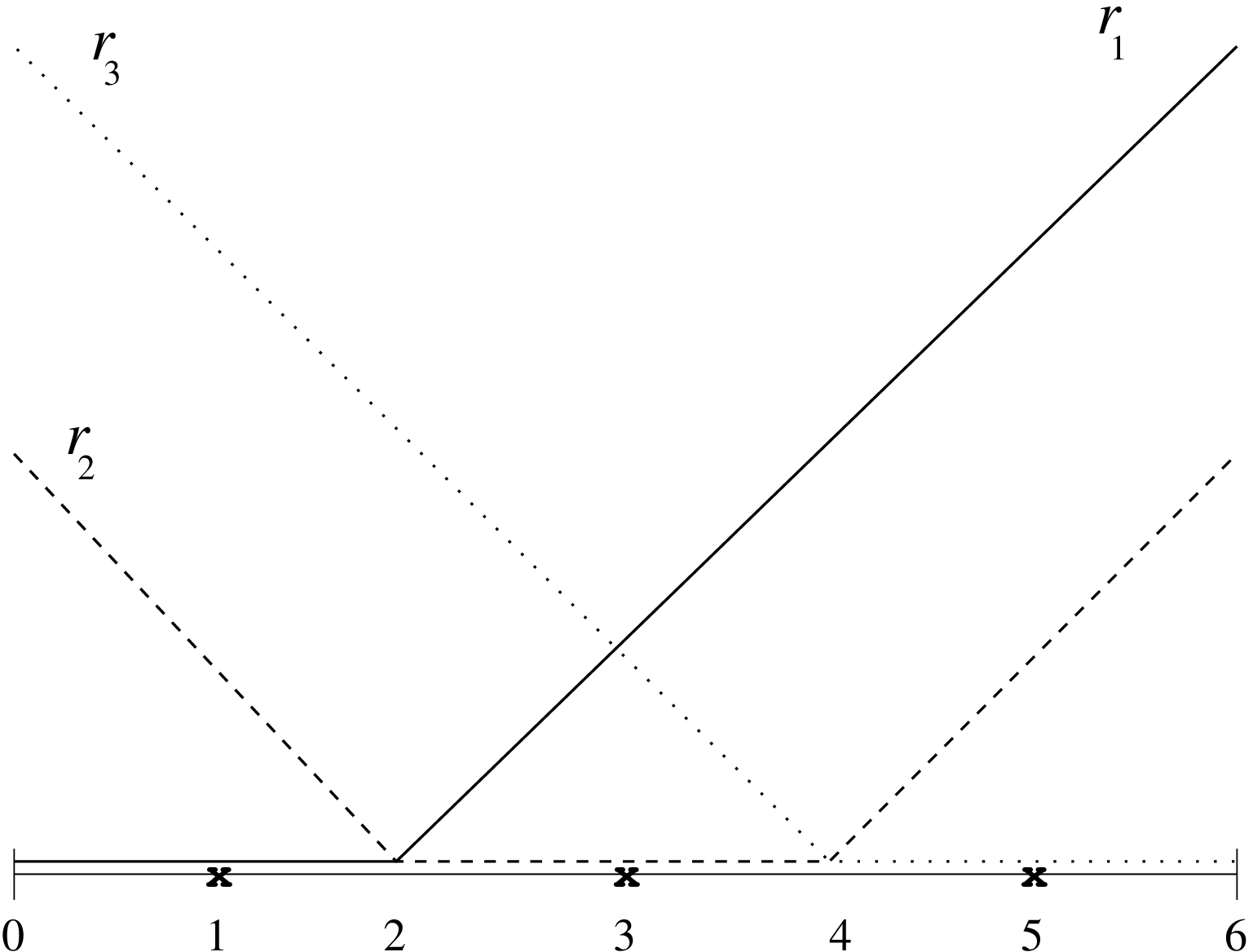}}
\caption{Distortion measure in Example~4.
Reproduction points are shown as {\bf x}'s.}
\end{figure}

Like in the discrete case, under some
additional assumptions on the rate-distortion
function $R(D)$, it is possible to get a 
stronger version of Theorem~3:

\vspace{0.1in}

{\sl Corollary 2.} Suppose $R(D)$ is 
differentiable and strictly convex 
on $(0,\Dmax)$. Under either one of 
the assumptions (a) and (b) in Theorem~3, 
there can be at most finitely many 
$D\in(0,\Dmax)$ such that $f(X_1)=0$ 
with probability one.

\vspace{0.1in}

{\sl Remark.} Under somewhat more restrictive
assumptions on the distortion measure $\rho$,
it is possible to prove that, for any $P$ with 
a continuous component as above, there can 
be at most $k(k+1)/2$ distortion levels $D$ for
which $f(X_1)=0$ with probability one. Since the
proof of this slightly stronger result relies 
on an argument different from the ones used to 
prove Theorems~2 and~3, we omit it here.

\section{Proofs}

\subsection{Preliminaries}

Before giving the proofs of Theorems~1, 2 and~3, 
we recall some definitions and notation from 
\cite{kontoyiannis-red:00} and give the 
precise form of the function $f$ (see equation 
(\ref{eq:eff}) below).

Let $P$ be a source distribution on $A$,
and let $Q$ be an arbitrary 
probability mass function on $\Ahat$.
Write $X$ for a random variable
with distribution $P$ on $A$,
and $Y$ for an independent random 
variable with distribution $Q$ on $\Ahat$.
Let $S =\{a\in\Ahat\;:\;Q(a)>0\}$
be the support of $Q$ 
and define
\ben
&\DminP& =\; E_{P}\left[ \min_{a\in S} 
		\rho(X,a)\right]\\
&\DmaxP\; & =\; E_{P\!\times\!Q}
	\left[\,\rho(X,Y)\right].
\een
For $\la\leq 0$, let
$$\LA_{P,Q}(\la)\,=\,
	E_{P}\left[
	\log_e
	E_{Q}
    	\left(
       	e^{\lambda \rho(X,Y)} 
    	\right)
  	\right],$$
and for $D\geq 0$ write $\LA^*_{P,Q}$
for the Fenchel-Legendre 
transform of $\LA_{P,Q}$,
        $$\LA^*_{P,Q}(D)=\sup_{\lambda\leq 0}\,
        [\la D-\LA_{P,Q}(\la)].$$
We also define
$$R(P,Q,D) \;=\;
\inf_{(X,Z)} [I(X;Z) + H(Q_Z\|Q)]$$
where $H(R\|Q)=\sum_{a\in \hat{A}}R(a)\log[R(a)/Q(a)]$
denotes the relative entropy (in bits)
between $R$ and $Q$,
$Q_{Z}$ denotes the distribution of $Z$,
and the infimum is over all jointly distributed
random variables $(X,Z)$ with values
in $A\!\times\!\Ahat$ such that
$X$ has distribution $P$ and
$E[\rho(X,Z)]\leq D$. In view
of (\ref{eq:rd-def}), 
we clearly have
\be
R(D)=\inf_{\mbox{\scriptsize all}\;Q} R(P,Q,D).
\label{eq:q-star}
\ee
In Lemma~1 and Proposition~1 
below we summarize some useful 
properties of $\LA_{P,Q}$,
$\LA^*_{P,Q}$ and
$R(P,Q,D)$
(see Lemma~1 and Propositions~1 and~2 
in \cite{kontoyiannis-red:00}).

\vspace{0.1in}

\begin{lemma}$\;$\\
\indent
$(i)$
	$\LA_{P,Q}$ is infinitely differentiable
        on $(-\infty,0)$,
	and $\LA''_{P,Q}(\la)\geq 0$ for all $\la\leq 0$.\\
\indent
$(ii)$
	If 
	$D\in(\DminP,\DmaxP)$
	then there exists a unique $\la<0$ such that 
	$\LA'_{P,Q}(\la)=D$ and
        $\LA^*_{P,Q}(D)=\la D -\LA_{P,Q}(\la)$.
\end{lemma}

\vspace{0.1in}

\begin{proposition}$\;$\\
\indent
$(i)$
	For all $D\geq 0$,
	$$R(P,Q,D) =  \inf_{W}\,
	E_{P}\left[ H(W(\cdot|X)\|Q(\cdot))\right],$$
	where the infimum is over all probability 
	measures $W$ on $A\!\times\!\Ahat$ such that
	the $A$-marginal of $W$ equals $P$
	and $E_{W}[\rho(X,Y)]\leq D.$\\
\indent
$(ii)$
	For all $D\geq 0$, 
	$\;R(P,Q,D)=(\log e)\LA^*_{P,Q}(D).$\\
\indent
$(iii)$
	For $0<D<\Dmax$ we have $0<R(D)<\infty$, whereas for 
	$D\geq\Dmax$, $R(D)=0$.\\
\indent
$(iv)$
	For every $D\in(0,\Dmax)$
	there exists a $Q=Q^*$ on $\Ahat$
	achieving the infimum in (\ref{eq:q-star}),
	and $D\in(D_{\rm min}^{P,Q^*},D_{\rm max}^{P,Q^*})$.
\end{proposition}

\vspace{0.1in}

For any distribution $P$ on $A$
and any distortion level $D\in(0,\Dmax(P))$,
by Proposition~1 we can pick
a $Q^*$ achieving the infimum 
in (\ref{eq:q-star}) so that
$R(D)=R(P,Q^*,D)$ and also
$D\in(D_{\rm min}^{P,Q^*},D_{\rm max}^{P,Q^*})$,
so by Lemma~1 we can pick a $\la^*<0$
with
\be
\la^*D-\LA_{P,Q^*}(\la^*) 
   &=& \LA^*_{P,Q^*}(D)
	\nonumber\\
   &=& (\log_e 2)R(P,Q^*,D)
	\nonumber\\
   &=& (\log_e 2)R(D).
\label{eq:represent}
\ee
Note also that 
\be
\la^*\to-\infty
\;\;\;\;\mbox{as}\;\;\;\;
D\to 0
\label{eq:slope}
\ee
(see the Appendix for a short proof).
Finally we can define the function $f$,
for $x\in A$,
\be
f(x)
\bydef
	(\log e)
	\left[
	\la^*D
	-
		\log_e E_{Q^*}\left(
                e^{\lambda^*\rho(x,Y)}
                \right)
	\right]
	-R(D).
\label{eq:eff}
\ee
Since $E_P[f(X_1)]=0$,
$f(X_1)=0$ with 
probability one if and only if
\be
\sum_{j=1}^kQ^*(a_j)
e^{\lambda^*\rho(x,a_j)}
	=\mbox{Constant,}
\;\;\mbox{for}\;P\mbox{-almost all $x$.}
\label{eq:condition}
\ee
Next we give an useful interpretation
for the constant $\la^*$ in the representation 
of $R(D)$ in (\ref{eq:represent}):
If $R(D)$ is differentiable at $D$, 
then $\la^*$ is proportional to its 
slope at $D$; Lemma~2 is proved in the Appendix.

\vspace{0.1in}

\begin{lemma}
For any $D\in(0,\Dmax)$:\\
\indent
$(i)$ We have
$(\log_e 2)R(D)=\sup_{\lambda\leq 0}\,
        [\la D-\Gamma(\la)],$
        where $\Gamma(\la)=\sup_{Q}\LA_{P,Q}(\la)$.\\
\indent
$(ii)$ Let
	$\la^*$ be chosen as in
	(\ref{eq:represent}). If $R(\cdot)$ is
	differentiable at $D$, 
	then $\la^*=(\log_e 2)R'(D).$
\end{lemma}

\subsection{Proofs in the Discrete Case}

For the proof of Theorem~1 we will need the
following lemma. 
It easily follows from Theorem~3.7 in Chapter~2 
of \cite{csiszar:book}
(see the Appendix). Recall the notation
$P_i=P(a_i)$ and $\rho_{ij}=\rho(a_i,a_j)$.

\vspace{0.1in}

\begin{lemma} 
A probability mass function
$Q^*$ on $A$ achieves the infimum in
(\ref{eq:q-star}) if and only if 
there exists a $\la^*<0$ such that the following
all hold:

(a) $\;$ $\LA'_{P,Q^*}(\la^*)=D.$

(b) $\;$ If we define, for $a_i,a_j\in A$,
	$$W(a_i,a_j)=P_iQ^*(a_j)\frac{e^{\lambda^*\rho_{ij}}}
		{\sum_{j'}Q^*(a_{j'})e^{\lambda^*\rho_{ij'}}
	 	}$$
	\hspace{0.5in}
	then the second marginal of $W$ is $Q^*$.

(c) $\;$ If $Q^*(a_j)=0$ for some $j$, then
	$$\sum_{i}P_i\frac{e^{\lambda^*\rho_{ij}}}
		{\sum_{j'}Q^*(a_{j'})e^{\lambda^*\rho_{ij'}}}
		\leq 1.$$
\end{lemma}

\vspace{0.1in}

\begin{example}
Here we present a simple
example illustrating the fact that
it may happen that $f(X_1)=0$ for a few
isolated values $D$ even when $P$ is
not uniform. Take $A=\Ahat=\{0,1,2\}$,
let $\alpha=\log_e[3e/(4-e)]$, and
consider the 
distortion measure
$$(\rho_{ij})
=
\left(
  \begin{array}{ccc}
  0 & 1 & \alpha\\
  1 & 0 & \alpha\\
  \alpha & \alpha & 0
  \end{array}
\right).
$$
Then with $P=Q^*=(4/13,4/13,5/13)$ and
$\la^*=-1$,
it is straightforward to check that 
condition~(b) of Lemma~3 holds (condition~(c) 
is irrelevant here),
and also (\ref{eq:condition}) is satisfied.
Therefore, at $D=\LA'_{P,Q^*}(\la^*)\approx 0.43$,
we must have $f(X_1)=0$ with probability one.
[Note, also, that the distortion measure
used here is not a permutation distortion
measure.]
\end{example}

\vspace{0.1in}

{\em Proof of Theorem 1, (a)}: Suppose 
$\rho$ is a permutation distortion measure
and $P$ is the uniform distribution on $A$, 
$P_i=1/k$ for all $i=1,\ldots,k$. First
we claim that for any $D\in(0,\Dmax)$
we can take $Q^*$ to also
be uniform. With $Q^*(a_j)=1/k$ for all $j$,
it suffices to find $\la^*<0$ satisfying
(a) and (b) of Lemma~3 (part (c) is 
irrelevant here).
We have
$D_{\rm min}^{P,Q^*}=0$
and
$$D_{\rm max}^{P,Q^*}=
\sum_{i,j}\frac{1}{k}\frac{1}{k}\rho_{ij}=
\frac{1}{k}\Sigma$$
where $\Sigma\bydef\sum_i\rho_{ij},$ 
which is independent of $j$ (since 
$\rho$ is a permutation). 
Also by the permutation property,
$\Dmax=\min_j E_P[\rho(X,a_j)]
=\min_j\sum_i(1/k)\rho_{ij}
=(1/k)\Sigma.$
Choose and fix a $D\in(0,\Dmax)$,
and 
pick $\la^*<0$ as in (\ref{eq:represent})
so that Lemma~3~(a) holds. With this $\la^*$ and
$Q^*$ being uniform let $W^*$ be as in Lemma~3~(b);
then
$$\sum_iW^*(a_i,a_j)
=\sum_i\frac{1}{k}\frac{1}{k}
	\frac{ e^{\lambda^*\rho_{ij}} }
	     { \sum_{j'}\frac{1}{k} e^{\lambda^*\rho_{ij'}} }
=\frac{1}{k}\sum_i
	\frac{ e^{\lambda^*\rho_{ij}} }
	     { \sum_{j'} e^{\lambda^*\rho_{ij'}} }.$$
But the sum in the denominator above 
\be
\sum_{j'} e^{\lambda^*\rho_{ij'}}
\;\;\;\;\mbox{is independent of $i$}
\label{eq:easy}
\ee
because $\rho$ is a 
permutation, so
$\sum_iW^*(a_i,a_j)=1/k=Q^*(a_j)$,
and (b) is satisfied. This proves that we can
take $Q^*$ to be uniform. Now
simply multiplying (\ref{eq:easy})
by $1/k$ we obtain (\ref{eq:condition}),
and this implies that $f(x)=0$ for all $x\in A$.
Since $D\in(0,\Dmax)$ was arbitrary, we are done.
\qed

\vspace{0.1in}

{\em Proof of Theorem 1, (b)}:
Let $D_n,$ $n\geq 1,$ be a sequence of
of distortion values in $(0,\Dmax)$
for which 
$f(X_1)=0$ with probability one,
and such that $D_n\downarrow 0$.
For each $D_n$, we can pick
$Q_n$ and a $\la_n<0$ 
as in (\ref{eq:represent})
such that 
$R(D_n)=(\log e)\LA_{P,Q_n}^*(\la_n)$
and $D_n=\LA'_{P,Q_n}(\la_n).$
Let 
$$\widetilde{D}=(\min_i P_i)(\min_{i\neq j} \rho_{ij})
>0.$$
Then for all $n$ large enough
so that $D_n<\widetilde{D}$, we must have 
$Q_n(a_i)>0$ for all $i$ (otherwise
it is trivial to check that 
$\LA'_{P,Q_n}(\la)\geq \widetilde{D}$
for any $\la<0$, contradicting the choice 
of $\la_n$). From now on we restrict 
attention to these large enough $n$'s.
As discussed above,
$f(X_1)=0$ with probability one
if and only if condition (\ref{eq:condition})
holds, which, in this case, becomes
\be
\sum_{j=1}^kQ_n(a_j)
e^{\lambda_n\rho_{ij}}
\;\;\;\;\mbox{is independent of $i$.}
\label{eq:condition1}
\ee
By Lemma~3~(b) we have that for all $j$
\ben
\sum_iP_i\frac{e^{\lambda_n\rho_{ij}}}
                {\sum_{j'}Q_n(a_{j'})e^{\lambda^*\rho_{ij'}}
                }=1,
\een
but by (\ref{eq:condition1})
the denominator is independent of $i$
so 
\be
\sum_iP_ie^{\lambda_n\rho_{ij}}=c_n,
\;\;\;\;\mbox{independent of $j$.}
\label{eq:simplified}
\ee
By (\ref{eq:slope}), $\la_n\to-\infty$ as $n\to\infty$,
so letting 
$n\to\infty$ yields
$P_j=\lim_n c_n$ for all $j,$
so $P$ is the uniform distribution
(recall our assumption that $\rho_{ij}=0$
if and only if $i=j$).
Moreover, from (\ref{eq:simplified})
it follows that
\be
\sum_ie^{\lambda_n\rho_{ij}} = kc_n,
\;\;\;\;\mbox{independent of $j$}.
\label{eq:independent}
\ee
To show that $\rho$ is a permutation,
fix two arbitrary indices $j\neq j'$ and
reorder the vectors 
$(\rho_{1j},\ldots,\rho_{kj})$
and $(\rho_{1j'},\ldots,\rho_{kj'})$
so that their elements are nondecreasing.
Let $(\sigma_{1},\ldots,\sigma_{k})$
and $(\sigma'_{1},\ldots,\sigma'_{k})$
be the 
corresponding ordered vectors.
Then
$\sigma_1=\sigma'_1=0$ and
(\ref{eq:independent}) implies that
$$\sum_{i=2}^ke^{\lambda_n(\sigma_i-\sigma_2')} =
\sum_{i=2}^ke^{\lambda_n(\sigma'_i-\sigma_2')}.$$
Next we show that if $\sigma_2\neq\sigma'_2$, say 
$\sigma_2>\sigma'_2$, we get a contradiction.
Since $\sigma_i-\sigma_2'>0$ for all $i\geq 2$,
the left-hand-side above tends to 0 
as $n\to\infty$,
but the right-hand-side
is $\geq 1$. Therefore
$\sigma_2=\sigma'_2$.
Continuing inductively, 
$\sigma_i=\sigma'_i$
for all $i$, so 
$(\rho_{1j},\ldots,\rho_{kj})$
and $(\rho_{1j'},\ldots,\rho_{kj'})$
are permutations of one another. Since
$j$ and $j'$ were arbitrary, this completes the proof.
\qed

\vspace{0.1in}

{\em Proof of Corollary 1}:
As before, let $D_n,$ $n\geq 1,$ 
be a sequence of of distortion 
values in $(0,\Dmax)$ for which
$f(X_1)=0$ with probability one,
and let
$Q_n$ and 
$\la_n<0$ be chosen such that
$R(D_n)=(\log e)\LA_{P,Q_n}^*(\la_n)$.
Since $R(D)$ is differentiable 
on $(0,\Dmax)$
(see \cite[Theorem~2.5.1]{berger:book}),
from Lemma~2 we get that 
$\la_n=(\log_e 2)R'(D_n)$.
Moreover,
since we assume that
$R(D)$ is strictly convex
on $(0,\Dmax)$, the $\la_n$ 
are all distinct. 

If the sequence $\{\la_n\}$ 
is unbounded, i.e., it has
a subsequence that
tends to $-\infty$, 
then we can proceed exactly 
as in the proof of Theorem~1.
So assume that the sequence 
$\{\la_n\}$ is bounded.
Since for each $n$, $R(P,Q_n,D_n)=R(D_n)>0$,
there must be a subset
$S$ of $\{1,2,\ldots,k\}$
of size $N$, say, with
$N=|S|\geq 2$, 
such that infinitely many of the $Q_n$
are supported on $\{a_j\,:\,j\in S\}$.
Without loss of generality we can relabel
the elements of $A$ so that 
$S=\{1,2,\ldots, N\}$.
If $N=k$ then we can again
repeat the argument
in the proof of Theorem~1.

Assuming $N\leq k-1,$ we proceed
to get a contradiction.
Since $f(x)=0$ with probability one,
condition (\ref{eq:condition}) implies
that 
$$
\sum_{j=1}^kQ_n(a_j)e^{\lambda_n\rho_{ij}}
=\sum_{j=1}^NQ_n(a_j)e^{\lambda_n\rho_{ij}}
=c_n,\;\;\mbox{for {\em all} $i$.}$$

Defining $\rho_{i0}=0$ for all $i$,
and letting $T(\la)$ denote the 
$(N+1)\!\times\!(N+1)$ matrix with 
entries $\exp(\lambda\rho_{ij})$
for $i=1,2,\ldots,N+1$ and $j=0,1,\ldots,N,$
the above conditions imply that 
$$T(\lambda_n)
\left(
  \begin{array}{c}
  -c_n\\
  Q_n(a_1)\\
  \vdots\\
  \vdots\\
  Q_n(a_N)
  \end{array}
\right)\;=\;
{\bf 0}\in\RL^{N+1}.
$$
Therefore $\det(T(\la_n))=0$
for all $\la_n$. The sequence
$\{\la_n\}$ is bounded so it must 
have an accumulation point, and
since $\det(T(\la))$ is an analytic 
function of $\la$ it can
only have isolated zeroes
unless it is identically zero
(see, e.g., the discussion in
\cite[Section~4.3.2]{ahlfors:book}).
So here we must have that $\det(T(\la))\equiv 0$ for
all $\la\leq 0$. But 
as $\la\to-\infty$, 
$T(\la)$
converges to the matrix
$$
T_\infty=
\left(
  \begin{array}{cccc}
  1 	 &   &           & \\
  1	 &   & {\bf I}_N & \\
  \vdots &   & 		 & \\
  1 	 & 0 & \cdots    & 0
  \end{array}
\right)
$$
which has determinant equal to 1 or $-1$
(${\bf I}_N$ denotes the $N\!\times\!N$ 
identity matrix),
and this provides the desired
contradiction.
\qed

\subsection{Proofs in the Continuous Case}

{\em Proof of Theorem 2}:
We argue by contradiction. Suppose
$f(X_1)=0$ with probability one for
some $D\in(0,\Dmax)$. Choose a
$Q^*$ and a $\la^*<0$ as in
(\ref{eq:represent}). 
Then (\ref{eq:condition}) implies that
$$
\sum_{j=1}^kQ^*(a_j)
e^{\lambda^* r_j(x)}
        =\mbox{Constant,}
\;\;\;\;\mbox{for}\;P\mbox{-almost all $x$,}$$
but since $P$ has an absolutely continuous
component with positive density on $I$,
and since the functions $r_j(\cdot)$ are
assumed to be continuous,
this holds for all $x\in I$,
and therefore contradicts the 
linear independence
assumption of Theorem~2.
\qed

\vspace{0.1in}

{\em Proof of Theorem 3}:
First we observe that condition (a)
immediately implies condition (b).
Therefore it suffices to show that 
if condition (b) holds, 
$f(X_1)$ cannot be equal
to zero with probability one for distortion
levels $D>0$ arbitrarily close to zero.
We proceed as in the proof of Corollary~1.
Assuming that there is a sequence
$D_n,$ $n\geq 1,$ 
of distortion values 
in $(0,\Dmax)$ for which
$f(X_1)=0$ with probability one,
and such that $D_n\downarrow 0$,
we will derive a contradiction.
Pick $Q_n$ and 
$\la_n<0$ such that
$R(D_n)=(\log e)\LA_{P,Q_n}^*(\la_n)$.
By (\ref{eq:condition}),
\be
\sum_{j=1}^kQ_n(a_j)
e^{\lambda_n r_j(x)} = c_n,
\;\;\;\;\mbox{for $P$-almost all $x\in I.$}
\label{eq:condition2}
\ee
Since 
$P$ has an absolutely continuous
component with positive density on $I$,
and since the functions $r_j(\cdot)$ are
assumed to be continuous,
(\ref{eq:condition2}) holds for all $x\in I$.
In particular, for the points $x_0,\ldots,x_k$
in condition (b), 
(\ref{eq:condition2}) becomes
$$\widetilde{T}(\la_n)\,(-c_n,Q_n(a_1),\ldots,Q_n(a_k))'=
{\bf 0}\in\RL^{k+1},$$
where $\widetilde{T}(\lambda)$ is the $(k+1)\!\times\!(k+1)$
matrix with entries $\exp(\la r_j(x_i))$,
$0\leq i,j\leq k$,
and $v'$ denotes the transpose of a vector $v$.
Therefore, since the entries of the
vector $(Q_n(a_1),\ldots,Q_n(a_k))$
sum to 1, it follows that 
$\det(\widetilde{T}(\la_n))=0$ for all $n$,
or, equivalently,
\be
\det(\widetilde{T}(\la_n))
&=&
\sum_{\pi}(-1)^{\signs(\pi)}
	e^{\lambda_n\sum_{j=0}^k r_j(x_{\pi(j)})}
	\nonumber\\
&=&
\sum_{\pi}(-1)^{\signs(\pi)}
	e^{\lambda_ns_{\pi}}=0,
\label{eq:determinant}
\ee
where the sums are taken over all permutations $\pi$ of the
set $\{0,1,\ldots,k\}$, and the constants $s_\pi$
are given by $\sum_{j=0}^k r_j(x_{\pi(j)})$.
Therefore, for any real number $s\geq 0$,
we must have that 
\be
\sum_{\pi\,:\,s_{\pi}=s}(-1)^{\signs(\pi)}=0.
\label{eq:contradiction}
\ee
To see this, let $\{s(1),s(2),\ldots\}$ be the (finite)
increasing sequence
of all possible values for the constants $s_{\pi}$.
Then (\ref{eq:determinant}) implies that 
$$\sum_{\pi\,:\,s_\pi=s(1)}(-1)^{\signs(\pi)}
	e^{\lambda_ns(1)}
   +\sum_{\pi\,:\,s_\pi>s(1)}(-1)^{\signs(\pi)}
	e^{\lambda_ns_{\pi}}=0.$$
By (\ref{eq:slope}), $\la_n\to-\infty$ as $n\to\infty$, 
so multiplying both sides by $e^{-\lambda_ns(1)}$
and letting $n\to\infty$ yields (\ref{eq:contradiction}) 
with $s=s(1)$. Continuing this way with $s(2)$,
then $s(3)$ and so on, proves (\ref{eq:contradiction})
for all $s$.

But now notice that condition (b) implies that,
if $\pi^*$ denotes the identity permutation,
then 
$s_{\pi}\neq s_{\pi^*}$ for all other 
permutations $\pi$. Therefore, taking
$s=s_{\pi^*}$ in (\ref{eq:contradiction})
we get the desired contradiction.\qed

\vspace{0.1in}

{\em Proof of Corollary 2}:
Let $D_n,$ $n\geq 1,$
be a sequence
of distortion values
in $(0,\Dmax)$ for which
$f(X_1)=0$ with probability one,
and pick $Q_n$ and
$\la_n<0$ as in the proof of Theorem~3.
If the sequence $\{\la_n\}$ is unbounded,
we can repeat the exact same proof
as for Theorem~3. So assume 
that $\{\la_n\}$ is bounded.
Since we also assume that $R(D)$ is
differentiable and strictly convex, it
follows from Lemma~2 that the 
$\la_n=(\log_e 2)R'(D_n)$ are all distinct.
Proceeding as in the proof of Theorem~3,
we get that
$\det(\widetilde{T}(\la))=0$
for all $\la=\la_n$.
The sequence
$\{\la_n\}$ is bounded so it must
have an accumulation point, and
$\det(\widetilde{T}(\la))$ is an analytic
function of $\la$. Therefore, arguing as in the
proof of Corollary~1,
$\det(T(\la))\equiv 0$ for
all $\la\leq 0$.
So we can find a sequence
$\la_m'\to-\infty$ for which 
$\det(\widetilde{T}(\la'_m))=0$.
With $\la'_m$ in place of $\la_n$,
the argument proceeds
exactly as in the proof of Theorem~3.
\qed

\section*{Acknowledgment}
We wish to thank the anonymous
Referees for their comments,
which helped improve the 
presentation 
of our paper.

\section*{Appendix}

{\em Proof of (\ref{eq:slope}):} Suppose
(\ref{eq:slope}) is false. Then it is
possible to
pick a constant $K<\infty$ and a sequence
of $D_n\in(0,\Dmax)$ with corresponding
$\la_n^*<0$, such that $D_n\to 0$ as $n\to\infty$
but $\la_n^*\geq -K$ for all $n$. 
Let $Q^*_n$
achieve (\ref{eq:q-star}) with $D=D_n$,
so that
\be
\LA'_{P,Q_n^*}(\la_n^*) = D_n.
\label{eq:lambda-star}
\ee
For each $n$,
recalling that $\rho(x,y)\leq M$
for all $x,y$,
\ben
\LA'_{P,Q_n}(\la_n^*)
&=&
   E_P\left[
        \frac{E_{Q_n}\left(\rho(X,Y)e^{\lambda^*_n\rho(X,Y)}\right)}
             {E_{Q_n}\left(e^{\lambda^*_n\rho(X,Y)}\right)}
        \right]\\
&\geq&
   E_P\left[E_{Q_n}\left(\rho(X,Y)e^{\lambda^*_n\rho(X,Y)}\right)\right]\\
&\geq&
   E_{Q_n}\left[ E_P\left(\rho(X,Y)e^{-KM}\right)\right]\\
&\geq&
   e^{-KM}\Dmax,
\een
which is bounded away from zero.
Since the $D_n\downarrow 0$, 
this contradicts (\ref{eq:lambda-star}).
\qed

\vspace{0.1in}

{\em Proof of Lemma 2:} Part $(i)$ immediately
follows from the minimax representation in 
\cite[Lemma~2]{kontoyiannis-red:00}.
For $(ii)$ note that, since
$\LA_{P,Q}(\la)$ is continuous and 
convex in $\la$ (Lemma~1), 
$\Gamma(\la)$ is lower semicontinuous 
and convex. Then by convex duality (see, e.g., 
Lemma 4.5.8 in \cite{dembo-zeitouni:book}),
it follows that $\Gamma(\la)=\sup_{x\geq 0}
[\la x- (\log_e 2)R(x)]$. For 
$D\in(0,\Dmax)$ and $\la^*$ as in 
(\ref{eq:represent}), we have
$$\Gamma(\la^*)=\la^*D-(\log_e 2)R(D)=\sup_{x\geq 0}
[\la^* x-(\log_e 2)R(x)].$$
But since 
$R(\cdot)$ is convex and (by assumption)
differentiable at $D$, it
must be that the derivative of
$[\la^* x-(\log_e 2)R(x)]$ vanishes
at $x=D$,
i.e., $\la^*=(\log_e 2)R'(D)$.
\qed

\vspace{0.1in}

{\em Proof of Lemma 3:} 
First suppose that for some $\la^*<0$, (a), (b) and (c) 
all hold. For $i=1,\ldots,k,$ let
$$B_i=\frac{P_i}
	{\sum_{j}Q^*(a_{j})e^{\lambda^*\rho_{ij}}}.$$
Then (b) and (c) imply that equations (3.19) and (3.20)
in \cite[p.~145]{csiszar:book} are satisfied with
$\delta=-\la^*$, so
by \cite[Theorem~3.7]{csiszar:book} equation (3.18)
is satisfied by $W^*$. This, together with Lemma~3.1 
in \cite[Chapter~2]{csiszar:book} implies that
$$R(D)=H(W\|P\!\times\!W^*_Y)$$
where $W^*_Y$ is
the second marginal of $W^*$. But 
$W_Y^*=Q^*$,
so $R(D)=E_P[H(W^*(\cdot|X)\|Q^*(\cdot))]$,
and 
by the definition of $W^*$ and Proposition~1,
$E_P[H(W^*(\cdot|X)\|Q^*(\cdot))]=R(P,Q^*,D)$.

Conversely, suppose $Q^*$ achieves the infimum
in (\ref{eq:q-star}). Then by Lemma~1 
there is a (unique) $\la^*<0$ such that (a)
holds, and letting $W^*$ be defined as in
(b) we also have
\ben
R(D)
&\eqa& R(P,Q^*,D)\\
&\eqb& H(W^*\|P\!\times\!Q^*)\\
&\eqc& H(W^*\|P\!\times\!W^*_Y) + H(W^*_Y\|Q^*)\\
&\geqd& H(W^*\|P\!\times\!W^*_Y)\\
&\geqe& R(D)
\een
where $(a)$ follows by assumption;
$(b)$ from (\ref{eq:represent}),
Proposition~1 and the
definition of $W^*$; $(c)$ by the
chain rule for relative entropy (see
\cite[Theorem~2.5.3]{cover:book});
$(d)$ is because relative entropy
is nonnegative; and $(e)$ follows
from the definition of $R(D)$ 
in (\ref{eq:rd-def}). Therefore
$H(W^*_Y\|Q^*)=0$, implying (b).
Finally note that 
the above argument shows 
that $W^*$ achieves $R(D)$.
Then by Theorem~3.7 in 
\cite[p.~145]{csiszar:book}
$W^*$ satisfies equation~(3.18) of
\cite[p.~145]{csiszar:book}
with $\delta=-\la^*$,
and by the uniqueness of 
the constants
$B_i$ and
equation~(3.19) 
of
\cite[p.~145]{csiszar:book}
we get (c).
\qed

\bibliographystyle{IEEE}

\begin{biography}{Ioannis Kontoyiannis}
was born in Athens, 
Greece, in 1972. He received the B.Sc. degree 
in mathematics in 1992 from Imperial College 
(University of London), and in 1993 he obtained 
a distinction in Part III of the Cambridge 
University Pure Mathematics Tripos. In 1997
he received the M.S. degree in statistics,
and in 1998 the Ph.D. degree in 
in electrical engineering,
both from Stanford University.
Between June and December 1995 he worked
at IBM Research, on a satellite image processing 
and compression project, funded by NASA 
and IBM. He has been with the Department 
of Statistics at Purdue University (and also, 
by courtesy, with the Department of Mathematics, 
and the School of Electrical and Computer 
Engineering) since 1998. 
During the 2000-01 academic year 
he is visiting the Applied Mathematics
Division of Brown University.
His research interests 
include data compression, 
applied probability,
statistical genetics,
nonparametric statistics,
entropy theory of stationary 
processes and random fields,
and ergodic theory.
\end{biography}

\begin{biography}{Amir Dembo}
received
the B.Sc. (Summa Cum Laude), and D.Sc. degrees in electrical engineering
from the Technion-Israel Institute of Technology, Haifa in 1980, 1986
respectively.  During 1980-1985 he was a Senior Research Engineer with
the Israel Defense Forces. In 1986/1987 he visited 
AT\&T Bell Laboratories and the Applied Mathematics 
Division of Brown University. Since 1988 he is
at Stanford University, where he is presently
an associate professor jointly in the Mathematics 
and Statistics departments. During 1994-1996 he 
was a professor of Electrical Engineering at the Technion, 
Israel.
He authored or co-authored more than 70 technical publications, 
including the book Large Deviations Techniques and Applications
(Second Edition, Springer-Verlag, 1998, with O. Zeitouni). 
From 1994-2000 he was an associate editor for the Annals of Probability,
and is currently serving on the editorial boards of the Annals of 
Applied Probability and of the Electronic Journal of Probability. 
He has worked in a number of areas including information and communication
theory, signal processing and estimation theory.
His current research interests are in probability theory and its 
applications to statistical physics, to queueing and communication systems, 
and to biomolecular data analysis.
\end{biography}

\end{document}